\documentclass[preprint,11pt,USletter]{elsarticle}
\usepackage{amsbsy}
\usepackage{amsfonts}
\usepackage{amsmath}
\usepackage{amssymb}
\usepackage[margin=1in]{geometry}

\newdefinition{dfn}{Definition}[section]
\newtheorem{theorem}{Theorem}[section]
\newtheorem{stw}{Proposition}[section]

\newproof{dwd}{Proof}
\newtheorem{wn}{Corollary}[section]
\newtheorem{lem}{Lemma}[section]
\newdefinition{uw}{Remark}[section]
\newdefinition{prob}{Problem}

\newdefinition{prk}{Example}

\UseRawInputEncoding

\begin{document}
\begin{frontmatter}
\title{On the bijective colouring of  Cantor trees based on  transducers}
\author[rvt]{Adam Woryna}
\ead{adam.woryna@polsl.pl}
\address[rvt]{Silesian University of Technology, Institute of Mathematics, ul. Kaszubska 23, 44-100 Gliwice, Poland}

\begin{abstract}
Given a vertex colouring  of the infinite $n$-ary Cantor tree with $m$ colours ($n,m\geq 2$), the natural problem arises: may this colouring  induce a bijective colouring of the infinite paths starting at the root, i.e., that  every  infinite  $m$-coloured string is used for some of these paths but different paths  are not coloured identically? In other words, we ask if the above vertex colouring may define  a bijective short map between the corresponding Cantor spaces. We show that the answer is  positive if and only if $n\geq m$, and provide an effective construction of the  bijective colouring  in terms of Mealy automata and functions defined by such automata. We also show that a finite Mealy automaton may  define such a  bijective colouring only in the trivial case, i.e. $m=n$.
\end{abstract}

\begin{keyword} Cantor tree \sep Bijection \sep Short Map \sep  Automaton \sep Function defined by automaton

\end{keyword}

\end{frontmatter}

\section{Introduction}

Let $X$ be a finite non-empty set ({\it alphabet}) with at least two elements and let  $X^*$ be the set of all finite words over $X$ together with the empty word $\epsilon$:
$$
X^*=\{\epsilon\}\cup\{x_1x_2\ldots x_l\colon x_i\in X,\;1\leq i\leq l,\;l=1,2,\ldots\}.
$$
Let  $X^\omega$ be the set of infinite words over $X$ (further  called  $\omega$-words):
$$
X^\omega=\{x_1x_2x_3\ldots\colon x_i\in X,\;i=1,2,\ldots\}.
$$
The  {\it $(n+1)$-ary} ($n\geq 1$) {\it Cantor tree} $T_X$ of  finite words over the alphabet  $X:=\{0, 1,\ldots, n\}$ is an example of an infinite, finite branching rooted tree, where the empty word $\epsilon$ is the root, and every vertex $w\in X^*$ has exactly $n+1$ children (immediate successors): $w0, w1, \ldots, wn$. The {\it boundary} $\partial T_X$ of the tree $T_X$ is the set of all infinite paths (without repetitions) that start at the root. These paths are in  one-to-one correspondence with  $\omega$-words over $X$.

We assume  that there are  infinitely many baubles each coloured with one of $m+1$  colours ($m\geq 1$) from the set  $Y:=\{0,1,\ldots, m\}$. Let us consider a decoration of the tree $T_X$ by  hanging  just one bauble on each vertex $w\neq \epsilon$. Such a  decoration naturally defines a colouring of every  infinite path from $\partial T_X$ with  some $\omega$-word over $Y$. We want to know if it is possible to decorate  $T_X$ in such a way that no two such paths have the same colouring and    every $\omega$-word over $Y$ is a colouring  of some of these paths. Let  $X^+=X^*\setminus\{\epsilon\}$ be the set of vertices  below the root. Then the problem can be formulated as follows.\\

\noindent
{\bf The colouring problem of the tree $T_X$:} is there a mapping $c\colon X^+\to Y$ (further called a {\it colouring} of the tree $T_X$) such that the map $\widetilde{c}\colon X^\omega\to Y^\omega$ defined as follows
$$
\widetilde{c}(x_1x_2x_3\ldots)=c(x_1)c(x_1x_2)c(x_1x_2x_3)\ldots,
$$
is a bijection between the sets  $X^\omega$ and $Y^\omega$?

If  $n=m$, then $X=Y=\{0, 1,\ldots, n\}$, and we can define a map $c\colon X^+\to Y$  as follows:
$$
c(w):=\mbox{\rm the last letter of }\;w,\;\;\;\;w\in X^+.
$$
Then  every  $\omega$-word $u=x_1x_2x_3\ldots \in X^\omega$ satisfies the equality:
$$
\widetilde{c}(u)=\widetilde{c}(x_1x_2x_3\ldots)=c(x_1)c(x_1x_2)c(x_1x_2x_3)\ldots=x_1x_2x_3\ldots=u.
$$
Hence $\widetilde{c}=Id_{X^\omega}$ and the above  formula  is an exemplary  solution of the problem. If  $n<m$, then  $X$ is a proper subset of  $Y$ and  for every mapping  $c\colon X^+\to Y$ the set $c(X)$  is also a proper subset of $Y$. In particular, the mapping  $\widetilde{c}\colon X^\omega\to Y^\omega$ cannot be onto, and we have a negative answer to the problem. The present paper shows that the problem has a positive answer in the non-trivial  case  $n>m$. We provide an  effective and neat construction of the required bijection $\widetilde{c}\colon X^\omega\to Y^\omega$ as a function defined by  a so-called finite asynchronous automaton, which   simultaneously can be  defined by an infinite Mealy automaton. The combinatorial language of such automata (also called transducers)  play an important  role in   various branches of mathematics. For example, in group theory, the self-similar transformations of the tree $T_X$ based on finite invertible  Mealy automata provide simple and elegant solutions of such  outstanding problems as Burnside's problem on periodic groups or Milnor's problem of group growth (\cite{6}). In the present paper, we also show that the required bijection $\widetilde{c}\colon X^\omega\to Y^\omega$ can be defined by a finite Mealy automaton  only in the trivial case, i.e. $m=n$.  It is also worth to note that the  Cantor spaces $X^\omega$ and $Y^\omega$ are homeomorphic according to the Brouwer's theorem~(\cite{2,7}). Thus our combinatorial solution determines when such a  homeomorphism between the spaces $X^\omega$ and $Y^\omega$ may be a 1-Lipschitz function, i.e. a short map.  For more on interesting properties of functions defined by  asynchronous and Mealy-type automata see
\cite{10,11,3,12,15,16,13,9,14}.

The paper is organized as follows. In Section~\ref{sec1}  the basic definitions concerning asynchronous and  Mealy-type automata are presented, as well as the concept of functions defined by these automata is recalled. Our description is similar to that of the survey paper \cite{6}.  In the proof of Proposition~\ref{stw1}, we  describe how colourings $c\colon X^+\to Y$ of the tree $T_X$ correspond to functions defined by Mealy automata. In Section~\ref{sec2}, we derive  (Proposition~\ref{g2}) the characterization of functions defined by Mealy automata as continuous functions $f\colon X^\omega\to Y^\omega$  which preserve the beginning of words, i.e. as short maps between the Cantor spaces $X^\omega$ and $Y^\omega$. As a  Corollary~\ref{cor22}, we obtain the condition when an asynchronous automaton can simulate a Mealy automaton, i.e. when a function defined by an asynchronous automaton can  simultaneously be defined by some Mealy automaton. In Section~\ref{sec3} the operation of composition of asynchronous automata is recalled, including  the proof that it corresponds to composition of mappings. The main construction and the solution of the problem is given in Section~\ref{sec4}, where we construct two families of asynchronous automata $A_n$ and $B_m$. In Propositions~\ref{ppp1}--\ref{kl1} we derive the formulae for transition and output functions in the composition of these automata. The main result is Theorem~\ref{twa3}, where we  show that if $n\geq m$, then some precisely constructed function defined by the composition of automata $A_n$ and $B_m$ is  bijective and can be defined by a Mealy automaton, i.e. it is a  bijective short map between the corresponding Cantor spaces. In the last section, we show (Thoerem~\ref{stst1} and Corollary~\ref{wnwn1}) that if $n\neq m$, then there is no  colouring $c\colon X^+\to Y$ such that the corresponding  function $\widetilde{c}\colon X^\omega\to Y^\omega$ is  bijective and may be  defined by a finite Mealy automaton.

\section{Automata and functions defined by automata}\label{sec1}

An {\it asynchronous automaton}  is a tuple
$$
A=(X,Q,Y,\varphi, \psi),
$$
where
\begin{itemize}
\item $X$ is an  {\it input alphabet},
\item $Q$ is a set of {\it internal states} of $A$,
\item $Y$ is an {\it output alphabet},
\item $\varphi\colon Q\times X\to Q$ is a {\it transition function} of $A$,
\item $\psi\colon Q\times X\to Y^*$ is an {\it output function} of $A$.
\end{itemize}
If the set $Q$ is finite, then we say that the automaton $A$ is {\it finite}. If the values of the output function  are just  one-letter words, that is $\psi$ is  of the form
$$
\psi\colon Q\times X\to Y,
$$
then the automaton $A$ is called a {\it Mealy automaton}.

It is convenient to present an arbitrary  automaton $A=(X, Q, Y, \varphi, \psi)$ by its Moore'a diagram, which is a directed and labeled  graph with the vertex set $Q$, and such that an arrow  goes from a vertex $q\in Q$  to a vertex $s\in Q$ if and only if $\varphi(q,x)=s$ for some letter $x\in X$. This arrow is labeled by a pair $x|v$, where $v=\psi(q,x)\in Y^*$ (the word $v$ may be empty). We then say that the automaton $A$ being in the state $q$ and reading the letter $x$ from the input tape, goes to the state $s=\varphi(q,x)$ and sends the word $v=\psi(q,x)$ on the output tape. A {\it path} in the automaton $A$ is a finite (perhaps empty) or infinite sequence of its arrows such that the end of each arrow from this sequence is the beginning of the next arrow (if  exists). If there is  a finite path
$$
q_0\stackrel{x_1|v_1}{\longrightarrow}q_1\stackrel{x_2|v_2}{\longrightarrow}\ldots q_{l-1}\stackrel{x_l|v_l}{\longrightarrow}q_l,
$$
then we say that  $A$ being in the initial state $q_0$ and reading the word $w=x_1x_2\ldots x_l\in X^*$ from the input tape, it sends the word $u=v_1v_2\ldots v_l\in Y^*$ on the output tape. Infinite paths  describe the action of  $A$ in the same way, that is when $A$ being in an initial state $q_0$ and  reading an infinite word $w=x_1x_2x_3\ldots\in X^\omega$ from the input tape, it  sends on the output tape a  word $u=v_1v_2v_3\ldots$ over the alphabet $Y$ (the word $u$  may be finite or infinite). Note that for every state $q\in Q$ and every word $w\in X^*\cup X^\omega$, the  automaton $A$ being in the  initial  state $q$ can read the word $w$ from the input tape. Then the unique  word $u\in Y^*\cup Y^\omega$  which is sent on the output tape, depends uniquely  on the  state $q$ and  on the word $w$. We denote this word by  $\overline{\psi}(q, w)$. For example, for the finite path as above, we have
$$
v_1v_2\ldots v_l=\overline{\psi}(q_0, x_1x_2\ldots x_l).
$$
If  $w\in X^*$ is finite, then  $\overline{\varphi}(q,w)\in Q$  denotes the state in which $A$ finishes its action after reading $w$ from the initial state $q$. In particular $\overline{\psi}(q,\epsilon)=\epsilon$ and $\overline{\varphi}(q,\epsilon)=q$ for all $q\in Q$. By the above action of  $A$, the following formulae should be clear:
\begin{eqnarray}
\overline{\varphi}(q, wv)&=&\overline{\varphi}(\overline{\varphi}(q,w), v),\label{f1}\\
\overline{\psi}(q, wu)&=&\overline{\psi}(q,w)\overline{\psi}(\overline{\varphi}(q,w), u)\label{f2}
\end{eqnarray}
for all $q\in Q$, $w,v\in X^*$, $u\in X^*\cup X^\omega$.

\begin{dfn}
The mapping
$$
f_q^A\colon X^*\cup X^\omega\to Y^*\cup Y^\omega,\;\;\;f_{q}^A(w)=\overline{\psi}(q, w)
$$
is called the {\it  function defined by the automaton $A$ in the state $q\in Q$}.
\end{dfn}

In particular, for all $w\in X^*$ and $u\in X^*\cup X^\omega$, we  have by~(\ref{f2}):
$$
f_{q}^A(wu)=f_q^A(w)f_{\overline{\varphi}(q, w)}^A(u).
$$
If $(w_i)_{i\geq 1}$ is an arbitrary infinite sequence  of finite words $w_i\in X^*$, then for their concatenation $w=w_1w_2w_3\ldots\in X^*\cup X^\omega$, we can write  by the formulae (\ref{f1})--(\ref{f2}):
\begin{equation}\label{f3}
f_q^A(w)=\overline{\psi}(q, w)=\overline{\psi}(q_0, w_1)\overline{\psi}(q_1, w_2)\overline{\psi}(q_2, w_3)\ldots,
\end{equation}
where
$$
q_0:=q,\;\;\;q_i:=\overline{\varphi}(q, w_1\ldots w_i)=\overline{\varphi}(q_{i-1}, w_{i}),\;i=1,2,\ldots.
$$

We denote by $f_{\omega, q}^A$ the restriction of  $f_{q}^A$ to the set  $X^\omega$. Obviously,  if $A$ is a Mealy automaton, then  $f_{\omega, q}^A$ sends every $\omega$-word over the alphabet $X$ into some $\omega$-word over the alphabet $Y$, that is we have in this case:
$$
f_{\omega, q}^A\colon X^\omega\to Y^\omega.
$$
In general, if $A=(X, Q, Y, \varphi, \psi)$ is an asynchronous automaton such that $f_{\omega, q}^A(u)\in Y^\omega$ for every $q\in Q$ and $u\in X^\omega$, then $A$ is called {\it nondegenerate}.

\begin{stw}\label{stw1}
For every Mealy automaton $A=(X, Q, Y, \varphi, \psi)$ and every state $q\in Q$ there is a colouring  $c\colon X^+\to Y$ of the tree $T_X$ such that $f_{\omega, q}^A=\widetilde{c}$. Conversely, for every colouring $c\colon X^+\to Y$ there is a Mealy automaton $A=(X, Q, Y, \varphi, \psi)$ and a state $q\in Q$ such that $f_{\omega, q}^A=\widetilde{c}$.
\end{stw}
\begin{dwd}
Let $A=(X, Q, Y, \varphi, \psi)$ be a Mealy automaton and $q\in Q$. Then for every $w\in X^*$ and $x\in X$ the word $\overline{\psi}(q,w)$ is a prefix of the word $\overline{\psi}(q,wx)$ and  the length  of  $\overline{\psi}(q,w)$ is the same as the length of $w$, i.e. $|\overline{\psi}(q,w)|=|w|$. Thus $\overline{\psi}(q,wx)-\overline{\psi}(q,w)\in Y$ (if a word $v$ is a prefix of a word $u$, then $u-v$ denotes the suffix of $u$ after removing $v$). Let us define  the colouring  $c\colon X^+\to Y$ as follows: if $w=x_1\ldots x_l\in X^+$ for some $x_i\in X$ ($1\leq i\leq l$, $l\geq 1$), then
$$
c(w)=c(x_1\ldots x_l)=\overline{\psi}(q,x_1\ldots x_l)-\overline{\psi}(q,x_1\ldots x_{l-1}).
$$
By an easy induction on $l\geq 1$, one can  show the equality:
$$
c(x_1\ldots x_l)=\psi(\overline{\varphi}(q, x_1\ldots x_{l-1}), x_l),\;\;\;l=1,2,3\ldots.
$$
Thus for every $\omega$-word $w=x_1x_2x_3\ldots\in X^\omega$, we have:
$$
\widetilde{c}(w)=c(x_1)c(x_1x_2)c(x_1x_2x_3)\ldots=\psi(q_0, x_1)\psi(q_1, x_2)\psi(q_2, x_3)\ldots,
$$
where $q_0:=q$, $q_i:=\overline{\varphi}(q, x_1\ldots x_i)$ for $i\geq 1$. Consequently $\widetilde{c}(w)=f_{\omega, q}^A(w)$.

Conversely, let $c\colon X^+\to Y$ be an arbitrary colouring of the tree $T_X$. Let us define a Mealy automaton $A=(X, Q, Y, \varphi, \psi)$ as follows:
$$
Q=X^*,\;\;\;\varphi(w,x)=wx,\;\;\;\psi(w,x)=c(wx)
$$
for all $w\in X^*$ and $x\in X$. Then for every $\omega$-word $w=x_1x_2x_3\ldots\in X^\omega$, we have:
\begin{eqnarray*}
\widetilde{c}(w)&=&c(x_1)c(x_1x_2)c(x_1x_2x_3)\ldots=\\
&=&\psi(\epsilon, x_1)\psi(x_1, x_2)\psi(x_1x_2,x_3)\ldots=\\
&=&\psi(q_0,x_1)\psi(q_1, x_2)\psi(q_2,x_3)\ldots,
\end{eqnarray*}
where $q_0=\epsilon$ and $q_i=x_1\ldots x_i$ for $i\geq 1$. Since
$q_i=q_{i-1}x_i=\varphi(q_{i-1}, x_i)$, by the fomula~(\ref{f3}), we obtain:
$\widetilde{c}(w)=f_{\omega, q_0}^A(x_1x_2x_3\ldots)=f_{\omega, q_0}^A(w)$.
\end{dwd}

\section{The characterization of functions $f_{\omega, q}^A$ defined by Mealy automata}\label{sec2}

For every $w\in X^*$ the {\it cone} $I_w\subseteq X^\omega$ corresponding to $w$ is the set  of all  infinite words having $w$ as a prefix:
$$
I_w:=wX^\omega=\{wu\colon u\in X^\omega\}=\{u\in X^\omega\colon w\preceq u\}.
$$
For all $w,w'\in X^*$, there are three possibilities:
$$
I_w\cap I_{w'}=\emptyset\;\;\;{\rm or}\;\;\;I_w\subseteq I_{w'}\;\;\;{\rm  or}\;\;\;I_{w'}\subseteq I_{w}.
$$
Moreover, the inclusion  $I_w\subseteq I_{w'}$ holds if and only if $w'$ is a prefix of $w$.

The set $X^\omega$ together with the cones $I_w$ ($w\in X^*$) form a topological space, which is homeomorphic to the classical Cantor set. This  topological space is induced by the metric $d_{X, \lambda}\colon X^\omega\times X^\omega\to\mathbb{R}$ defined for any real number $0<\lambda<1$ as follows:
$$
d_{X, \lambda}(u_1, u_2)=\lambda^n,
$$
where $n$ is the length of the longest common prefix of $u_1$ and $u_2$ (if $u_1=u_2$, then we assume $d_{X, \lambda}(u_1, u_2)=0$). One can show that $d_{X,\lambda}$ satisfies  the inequality
$$
d_{X, \lambda}(u_1, u_3)\leq \max\{d_{X, \lambda}(u_1, u_2), d_{X, \lambda}(u_2,u_3)\}
$$
for all $u_1, u_2, u_3\in X^\omega$ (\cite{7}), which means that $d_{X,\lambda}$ is also an ultrametric. The cones $I_w$ are open balls in the metric space $(X^\omega, d_{X,\lambda})$. According to the Brouwer's characterization of Cantor spaces, the space $X^\omega$ is the only (up to homeomorphism) perfect non-empty compact metrizable space which is Hausdorff and has a basis consisting of clopen sets (\cite{7}).

\begin{theorem}[\cite{6}, Theorem~2.4]
For any two alphabets $X$ and $Y$ and a function $f\colon X^\omega\to Y^\omega$ the following two statements are equivalent.
\begin{itemize}
\item[(i)] The function $f$ is continuous.
\item[(ii)] There is an asynchronous automaton $A=(X, Q, Y, \varphi, \psi)$ such that $f=f_{\omega, q}^A$ for some $q\in Q$.
\end{itemize}
\end{theorem}

Thus the class of functions  $f\colon X^\omega\to Y^\omega$ which are continuous corresponds to the class of functions defined by asynchronous automata. We distinguish in this  class the functions defined by Mealy automata, i.e. the  functions $f\colon X^\omega\to Y^\omega$ for which there is a Mealy automaton $A=(X, Q, Y, \varphi, \psi)$ such that $f=f_{\omega, q}^A$ for some  $q\in Q$. In particular, every function   $f\colon X^\omega\to Y^\omega$ defined by a Mealy automaton is continuous, but the converse is not true. For example, the {\it unilateral shift}  $f\colon X^\omega\to X^\omega$ defined as follows:
$$
f(x_1x_2x_3\ldots)=x_2x_3x_4\ldots,
$$
is continuous, but there is no Mealy automaton  defining $f$. Indeed, if this function would be defined by a Mealy automaton $A=(X, Q, X, \varphi, \psi)$, then by the formula (\ref{f3}),  for any two different letters $x, x'\in X$ both the words $f(x^\infty)$ and $f(xx'x^\infty)$ will start with the same letter (equal to $\psi(q,x)$ for some state $q\in Q$). But, according to the definition of $f$, the word $f(x^\infty)$ starts with $x$ and the word $f(xx'x^\infty)$ starts with $x'$.

For any subset $S\subseteq X^*\cup X^\omega$ let  $P(S)$ denotes the longest common prefix of the words in $S$. Note that the word  $P(S)$ is infinite if and only if  $S$ is a one-element set consisting of an $\omega$-word. Otherwise, we have  $P(S)\in X^*$.

\begin{dfn}
For an arbitrary function  $f\colon X^\omega\to Y^\omega$, we say that  $f$ {\it preserves the beginning of words}, if for all $u,v\in X^\omega$ the longest common prefix of the images  $f(u)$ and $f(v)$ is not shorter that the longest common prefix of $u$ and $v$, that is
$|P(f(u),f(v))|\geq |P(u,v)|$.
\end{dfn}

\begin{uw}
According to the above definition, a  function  $f\colon X^\omega\to Y^\omega$  preserves the beginning of words if and only if $d_{Y, \lambda}(f(u), f(v))\leq d_{X, \lambda}(u,v)$ for all  $u,v\in X^\omega$, which means that  $f$ is a short map between the spaces $(X^\omega, d_{X,\lambda})$ and $(Y^\omega, d_{Y,\lambda})$.
\end{uw}

Let  $f\colon X^\omega\to Y^\omega$ be an arbitrary function that preserves the beginning of words. Fix $w\in X^*$. Then for all $u,u'\in X^\omega$ the words $f(wu)$ and $f(wu')$ have common prefix of length $|w|$. In particular, the longest common prefix of the image
$f(I_w)=f(wX^\omega)\subseteq Y^\omega$
is not shorter than $w$, that is
$|P(f(I_w))|\geq |w|$. Thus there is a unique word $v\in Y^*$ for which $|v|=|w|$ and $v\preceq f(wu)$ for every  $u\in X^\omega$. In particular, $f(wu)\in I_v$ for every  $u\in X^\omega$. Thus the word  $v$ is the only word over $Y$ such that following two conditions hold (i) $|v|=|w|$, (ii) $f(I_w)\subseteq I_v$. Let us denote this word by $L_f(w)$.
Obviously, we have $L_f(\epsilon)=\epsilon$. The word  $L_f(w)$  is not necessarily  the longest common prefix of the words in the image $f(I_w)$, but it is  equal to the prefix of length $|w|$ of any word   $u\in f(I_w)$.

\begin{uw}
For all $x\in X$ and $w\in X^*$ the word  $L_f(w)$ is a prefix of the word  $L_f(wx)$. Indeed, if we  denote $v:=L_f(w)$ and  $v':=L_f(wx)$, then $v'$ is the unique word of length  $|v'|=|w|+1=|v|+1$ which satisfies: $f(I_{wx})\subseteq I_{v'}$. Since $I_{wx}\subseteq I_{w}$, we have: $f(I_{wx})\subseteq f(I_{w})\subseteq I_v$. Thus $I_v\cap I_{v'}\neq \emptyset$ and hence $I_{v'}\subseteq I_v$. Consequently $v\prec v'$.
\end{uw}

\begin{stw}\label{g2}
A function $f\colon X^\omega\to Y^\omega$ is defined by a Mealy automaton  if and only if it preserves the beginning of words, that is $f$ is a short map between the spaces $(X^\omega, d_{X, \lambda})$ and $(Y^\omega, d_{Y, \lambda})$.
\end{stw}
\begin{dwd}
If $f=f_{\omega, q}^A$ for some  Mealy automaton $A=(X, Q, Y, \varphi, \psi)$ and its state  $q\in Q$, then we see by  the formula (\ref{f3}) that $f$ preserves the beginning of words. Conversely, assume that a function $f\colon X^\omega\to Y^\omega$ preserves the beginning of words. For all $w\in X^*$ and $x\in X$ let
$L_f(wx)-L_f(w)\in Y$
be a one-letter suffix of the word $L_f(wx)$ when removing the prefix $L_f(w)$. Let us define a Mealy automaton
$A=(X, X^*, Y, \varphi, \psi)$ as follows:
$$
\varphi(w, x)=wx,\;\;\;\psi(w, x)=L_f(wx)-L_f(w).
$$
Fix an arbitrarily infinite word $x_1x_2x_3\ldots\in X^\omega$ and let us denote:
$$
f(x_1x_2x_3\ldots)=y_1y_2y_3\ldots\in Y^\omega,\;\;\;y_i\in Y.
$$
For every  $i\geq 1$ the word $L_f(x_1\ldots x_i)\in Y^*$ is a common prefix of length  $i$ of all  words in the image $f(I_{x_1\ldots x_i})$. Since  the word $f(x_1x_2\ldots)=y_1y_2y_3\ldots$ belongs to this image, we obtain:
$L_f(x_1\ldots x_i)=y_1\ldots y_i$ for every $i\geq 1$. By the formula (\ref{f3}), we can write:
$$
f^A_{\omega,\epsilon}(x_1x_2x_3\ldots)=\psi(q_0, x_1)\psi(q_1, x_2)\psi(q_2, x_3)\ldots,
$$
where $q_0:=\epsilon$, $q_i:=\overline{\varphi}(q_0, x_1\ldots x_{i})=x_1\ldots x_i$ for $i>0$. By the definition of the output function  $\psi$, we have:
$$
\psi(q_0, x_1)=\psi(\epsilon, x_1)=L_f(x_1)-L_f(\epsilon)=L_f(x_1)=y_1.
$$
For every  $i\geq 1$, we also have:
\begin{eqnarray*}
\psi(q_i, x_{i+1})&=&\psi(x_1\ldots x_i, x_{i+1})=\\
&=&L_f(x_1\ldots x_{i+1})-L_f(x_1\ldots x_i)=y_1\ldots y_{i+1}-y_1\ldots y_i=y_{i+1}.
\end{eqnarray*}
Thus
$$
f^A_{\omega,\epsilon}(x_1x_2x_3\ldots)=\psi(q_0, x_1)\psi(q_1, x_2)\psi(q_2, x_3)\ldots=y_1y_2y_3\ldots=f(x_1x_2x_3\ldots).
$$
In consequence $f=f_{\omega,\epsilon}^A$, which means that $f$ is defined by a Mealy automaton $A$.
\end{dwd}

By Proposition~\ref{stw1} and the proof of Proposition~\ref{g2}, we have:

\begin{wn}
The colouring problem of the tree $T_X$ has a positive solution if and only if there is a bijective function $f\colon X^\omega\to Y^\omega$ defined by a Mealy automaton. For such a function $f$, we have: $f=\widetilde{c}$, where the colouring $c\colon X^+\to Y$  is defined as follows: $c(wx)=L_f(wx)-L_f(w)$ for all $w\in X^*$ and $x\in X$.
\end{wn}

By Proposition~\ref{g2} and  formula~(\ref{f3}), we obtain:

\begin{wn}\label{cor22}
Let  $A=(X, Q, Y, \varphi, \psi)$ be an asynchronous automaton such that the word $\psi(q,x)\in Y^*$ is non-empty for all $q\in Q$ and $x\in X$. Then for every $q\in Q$ the function  $f=f_{\omega, q}^A\colon X^\omega\to Y^\omega$ is simultaneously defined by some  Mealy automaton. However, the colouring $c\colon X^+\to Y$ for which  $f=\widetilde{c}$, can be determined directly from the action of the asynchronous  automaton $A$ on  finite words in the following way: for any word $w=x_1\ldots x_n\in X^+$,  the letter $c(w)\in Y$ is equal to the $n$-th letter of the word $\overline{\psi}(q,w)$, i.e. to the $n$-th letter of the word  sent by $A$ on the output tape after reading  $w$ from the initial state $q$.
\end{wn}

\section{The composition of asynchronous automata}\label{sec3}

Let $A=(X, Q^A, Y, \varphi^A, \psi^A)$ and $B=(Y, Q^B, Z, \varphi^B, \psi^B)$ be arbitrary asynchronous automata. The automaton $C=(X, Q^A\times Q^B, Z, \varphi, \psi)$ defined as follows
\begin{eqnarray*}
\varphi((q,s),x)&=&(\varphi^A(q,x),\overline{\varphi}^B(s,\psi^A(q,x))),\\
\psi((q,s),x)&=&\overline{\psi}^B(s,\psi^A(q,x))
\end{eqnarray*}
for all $(q,s)\in Q^A\times Q^B$, $x\in X$ is called the {\it composition of automata $A$ and $B$} and is denoted by $A\circ B$.

\begin{uw}
In the below proposition and further, we use the right action convention for composition of mappings,  that is if $f\colon X\to Y$ and $g\colon Y\to Z$, then  $f\circ g(x)=g(f(x))$ for every $x\in X$.
\end{uw}

\begin{stw}\label{ppp2}
Let $A=(X, Q^A, Y, \varphi^A, \psi^A)$ and $B=(Y, Q^B, Z, \varphi^B, \psi^B)$ be  asynchronous automata, and the automaton  $A$ is nondegenerate.  Then  for every $(q,s)\in Q^A\times Q^B$, we have:
$$
f_{\omega, (q,s)}^{A\circ B}=f_{\omega, q}^A\circ f_{\omega, s}^B.
$$
\end{stw}
\begin{dwd}
Let $u=x_1x_2x_3\ldots\in X^\omega$ be arbitrary. We need to show the equality $f_{\omega, (q,s)}^{A\circ B}(u)=f_{\omega, q}^A\circ f_{\omega, s}^B(u)$. We have
\begin{eqnarray*}
f_{\omega, q}^A\circ f_{\omega, s}^B(u)=f_{\omega, s}^B(f_{\omega, q}^A(u))=f_{\omega, s}^B\left(\overline{\psi}^A(q,u)\right)=\\
=\overline{\psi}^B\left(s,\overline{\psi}^A(q,u)\right)=\overline{\psi}^B\left(s,\psi^A(q_0, x_1)\psi^A(q_1, x_2)\ldots\right),
\end{eqnarray*}
where
$$
q_0=q,\;\;\;q_i=\varphi^A(q_{i-1},x_i),\;\;\;i=1,2,\ldots.
$$
If we denote:
\begin{equation}\label{ee0}
y_i=\psi^A(q_{i-1}, x_i),\;\;\;i=1,2,\ldots,
\end{equation}
then  $y_i\in Y^*$ and we can write:
\begin{equation}\label{ee1}
f_{\omega, q}^A\circ f_{\omega, s}^B(u)=\overline{\psi}^B(s,y_1y_2\ldots)=\overline{\psi}^B(s_0, y_1)\overline{\psi}^B(s_1, y_2)\ldots,
\end{equation}
where
$$
s_0=s,\;\;\;s_i=\overline{\varphi}^B(s_{i-1},y_i),\;\;\;i=1,2,\ldots.
$$
On the other hand, we have by the definition of the composition $A\circ B$:
\begin{equation}\label{ee2}
f_{\omega, (q,s)}^{A\circ B}(u)=\overline{\psi}((q,s),x_1x_2\ldots)=\psi((t_0, r_0), x_1)\psi((t_1, r_1), x_1)\ldots,
\end{equation}
where
$$
(t_0, r_0)=(q,s),\;\;\;(t_i, r_i)=\varphi((t_{i-1}, r_{i-1}), x_i),\;\;\;i=1,2,\ldots.
$$
By the definition of the mapping $\varphi$, we have for every $i\geq 1$:
$$
(t_i, r_i)=\left(\varphi^A(t_{i-1}, x_i), \overline{\varphi}^B(r_{i-1}, \psi^A(t_{i-1}, x_i))\right).
$$
In particular $t_0=q$ and $t_{i}=\varphi^A(t_{i-1}, x_i)$ for $i=1,2\ldots.$ Consequently
\begin{equation}\label{ee3}
t_i=q_i,\;\;\;i=0,1,\ldots.
\end{equation}
We also have: $r_0=s$ and for every $i\geq 1$:
$$
r_i=\overline{\varphi}^B(r_{i-1}, \psi^A(t_{i-1}, x_i))=\overline{\varphi}^B(r_{i-1}, \psi^A(q_{i-1}, x_i))=\overline{\varphi}^B(r_{i-1}, y_i),
$$
which implies:
\begin{equation}\label{ee4}
r_i=s_i,\;\;\;i=0,1,\ldots.
\end{equation}
Now, by (\ref{ee1})--(\ref{ee2}) and by (\ref{ee3})--(\ref{ee4}), it is enough to show the equalities:
$$
\psi((q_{i-1}, s_{i-1}), x_i)=\overline{\psi}^B(s_{i-1}, y_i),\;\;\;i=1,2,\ldots.
$$
But these equalities  directly follow from (\ref{ee0}) and the definition of $\psi$.
\end{dwd}

\section{The main construction: automata $A_n$, $B_m$ and their composition}\label{sec4}

For every $n\in\{1,2,\ldots\}$ let
$$
A_n=(X, \{\sigma\}, Y, \varphi^{A_n}, \psi^{A_n})
$$
be the asynchronous automaton with a single state $\sigma$, the input alphabet $X=\{0,1,\ldots, n\}$, the output alphabet $Y=\{0,1\}$ and the output function $\psi^{A_n}\colon \{\sigma\}\times \{0,1,\ldots, n\}\to \{0,1\}^*$ defined as follows:
$$
\psi^{A_n}(\sigma, x)=\left\{
\begin{array}{ll}
0^x1,&0\leq x\leq n-1,\\
0^n,&x=n.
\end{array}
\right.
$$
\begin{figure}[hbtp]
\centering
\includegraphics[width=5cm]{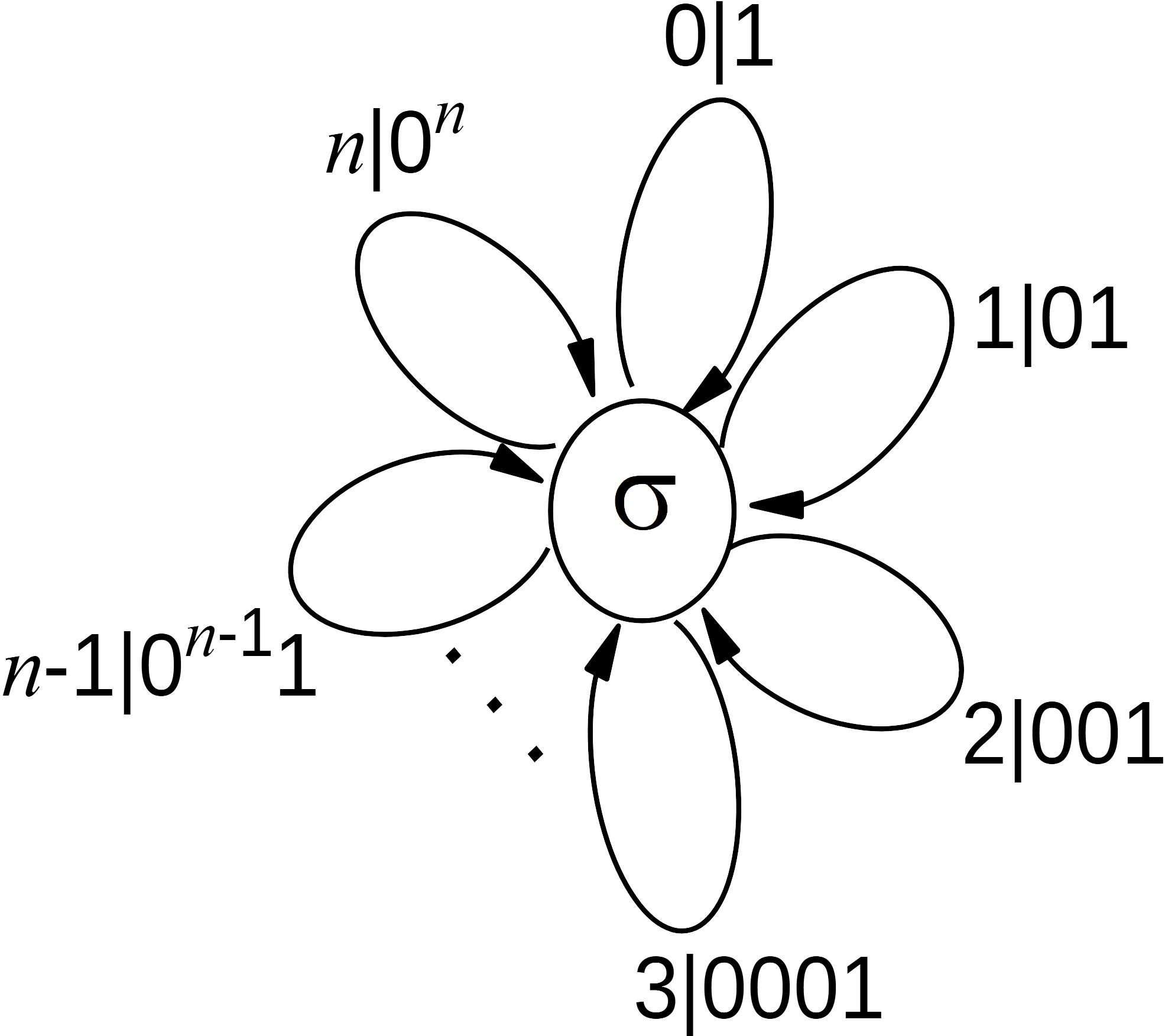}
\caption{The automaton $A_n$}
\end{figure}

\begin{figure}[hbtp]
\centering
\includegraphics[width=12cm]{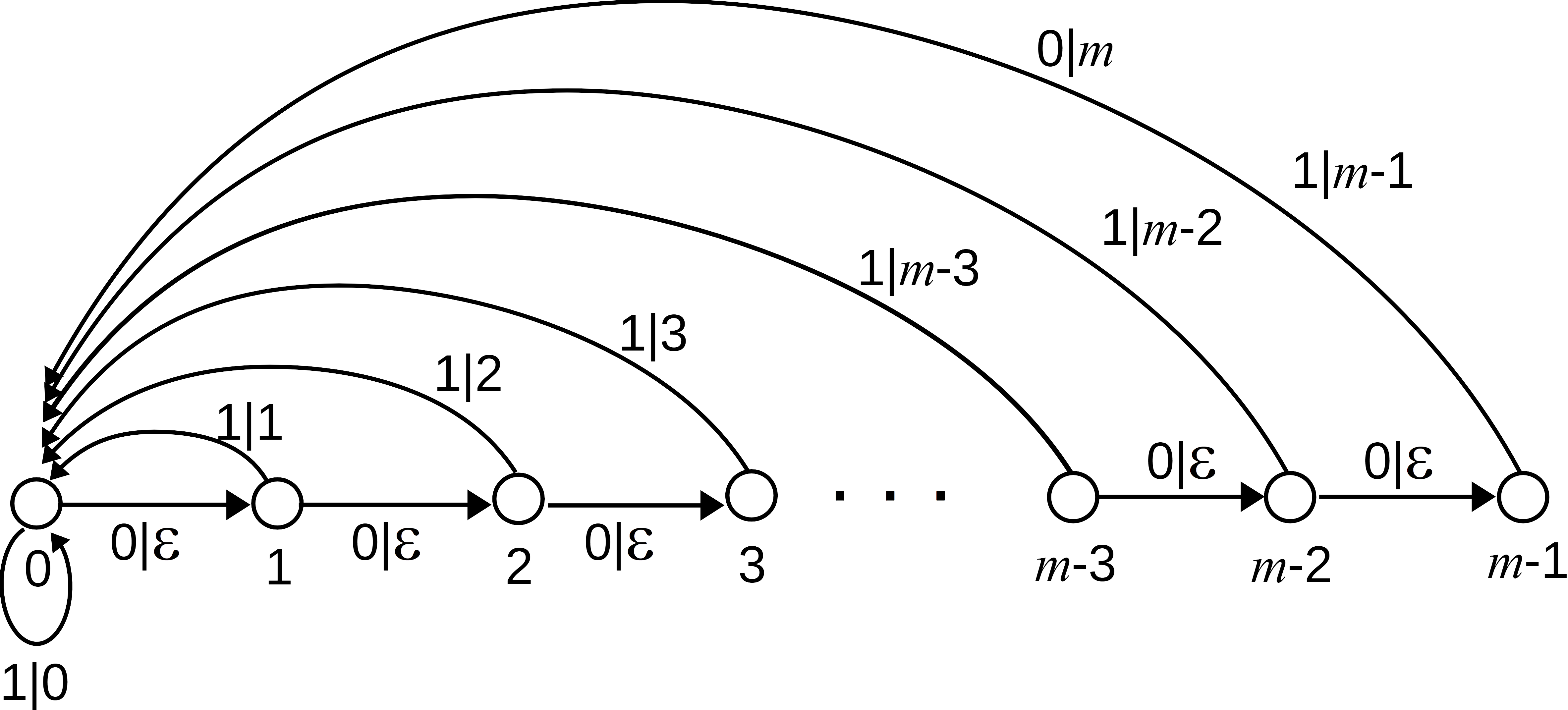}
\caption{The automaton $B_m$}
\label{fig2}
\end{figure}

For every $m\in\{1,2,\ldots\}$ we also define the automaton
$$
B_m=(Y, Q^{B_m}, Z, \varphi^{B_m}, \psi^{B_m})
$$
with the input alphabet $Y=\{0,1\}$, the output alphabet $Z=\{0,1,\ldots, m\}$, the set of states
$$
Q^{B_m}=\{0,1,\ldots, m-1\},
$$
and the transition and output functions $\varphi^{B_m}$, $\psi^{B_m}$ defined as follows:
\begin{eqnarray*}
\varphi^{B_m}(q, x)&=&\left\{
\begin{array}{ll}
0,&x=1,\\
q+_m1,&x=0,
\end{array}
\right.\\
\psi^{B_m}(q, x)&=&\left\{
\begin{array}{ll}
q,&x=1,\\
\epsilon,&x=0\;\;{\rm and}\;\;q\neq m-1,\\
m,&x=0\;\;{\rm and}\;\;q=m-1,
\end{array}
\right.
\end{eqnarray*}
where $+_m$ denotes the addition (mod $m$). The Moore'a diagram of the automaton $B_m$ is depicted in Figure~\ref{fig2}.

Note that  the automata $A_n$ and $B_m$ are nondegenerate. It is clear for the automaton $A_n$, as the values of its output function $\psi^{A_n}$ are non-empty. The automaton $B_m$ is also nondegenerate. Indeed, for every state $q\in Q^{B_m}$ the word $\psi^B(q,1)$ is nonempty (because   equal $q$). This implies that  if a word $u\in \{0,1\}^\omega$ contains infinitely many $1$'s, then for every state $q\in Q^{B_m}$ the word $f_{\omega, q}^{B_m}(u)$ is infinite. If $u=w0^\infty$ for some $w\in \{0,1\}^*$, then we have
$$
f_{\omega, q}^{B_m}(u)=\overline{\psi}^{B_m}(q,w)\overline{\psi}^{B_m}(q', 0^\infty),
$$
where $q'=\overline{\varphi}^{B_m}(q, w)$. But, directly from the Moore'a diagram  of  $B_m$, we see that if $t\geq m$, then for every state $s\in Q^{B_m}$ the word
$\overline{\psi}^{B_m}(s, 0^t)$ contains at least one letter $m$, because the corresponding path with the initial state $s$ contains the arrow
$$
m-1\stackrel{0|m}{\longrightarrow}0.
$$
In particular, the word $\overline{\psi}^{B_m}(q', 0^\infty)$ is infinite, and hence the word $f_{\omega, q}^{B_m}(u)$ is also infinite.

\begin{lem}\label{lem1}
For every $t\geq 0$ and $q\in Q^{B_m}$, we have:
\begin{equation}\label{pf1}
\overline{\varphi}^{B_m}(q, 0^t)=q+_mt,\;\;\;\overline{\varphi}^{B_m}(q, 0^t1)=0,
\end{equation}
and
\begin{equation}\label{pf2}
\overline{\psi}^{B_m}(q, 0^t)=m^{\lfloor\frac{q+t}{m}\rfloor},\;\;\;\overline{\psi}^{B_m}(q, 0^t1)=m^{\lfloor\frac{q+t}{m}\rfloor}(q+_mt).
\end{equation}
\end{lem}
\begin{dwd}
The formulae~(\ref{pf1}) follows directly from the definition of $\varphi^{B_m}$ and from the observation that $\varphi^{B_m}(q,1)=0$ for every state $q\in Q^B$. To show~(\ref{pf2}), we have:
$$
\overline{\psi}^{B_m}(q,0^t)=z_1z_2\ldots z_t,\;\;\;z_i=\psi^{B_m}(q_i, 0),\;\;\;i=1,2,\ldots, t,
$$
where
$$
q_i=\overline{\varphi}^{B_m}(q, 0^{i-1})=q+_m(i-1).
$$
Thus, by the definition of $\psi^{B_m}$, we have: $z_i\in\{\epsilon, m\}$, and $z_i=m$  if and only if $q+_m(i-1)=m-1$. Thus
$$
\overline{\psi}^{B_m}(q,0^t)=m^r,
$$
where $r\geq 0$ is the number of those elements from the sequence
$$
q,\;q+_m1,\;q+_m2,\;\ldots, q+_m(t-1)
$$
which are equal to $m-1$. Since $q\in \{0,1,\ldots, m-1\}$, these elements form the following subsequence (perhaps empty):
$$
q+_m(m-1-q),\;q+_m(2m-1-q),\;\ldots,\;q+_m(km-1-q),
$$
where the number $k\geq 0$ satisfies:
$$
km-1-q\leq t-1,\;\;\;(k+1)m-1-q>t-1.
$$
The  two last inequalities imply $k=\lfloor (q+t)/m\rfloor$. Since the above subsequence is of length $k$, we obtain:
$$
r=k=\lfloor (q+t)/m\rfloor.
$$
Now, since $\psi^{B_m}(s,1)=s$ for every $s\in Q^{B_m}$, we can write
\begin{eqnarray*}
\overline{\psi}^{B_m}(q,0^t1)=\overline{\psi}^{B_m}(q,0^t)\overline{\psi}^{B_m}(\overline{\varphi}^{B_m}(q,0^t),1)=\\
=m^{\lfloor (q+t)/m\rfloor}\overline{\psi}^{B_m}(q+_mt,1)=m^{\lfloor (q+t)/m\rfloor}(q+_mt).
\end{eqnarray*}
\end{dwd}

\begin{stw}\label{ppp1}
For all $n,m\geq 1$,  the composition  $A_n\circ B_m=(X, Q, Z, \varphi, \psi)$ of automata $A_n$ and $B_m$ is defined as follows:
\begin{itemize}
\item $X=\{0,1,\ldots, n\}$,
\item $Q=\{(\sigma,0),\ldots, (\sigma, m-1)\}$,
\item $Z=\{0,1,\ldots, m\}$,
\item $\varphi((\sigma, q), x)=\left\{
\begin{array}{ll}
(\sigma, 0),&0\leq x\leq n-1,\\
(\sigma, q+_mn),&x=n,
\end{array}\right.$
\item $\psi((\sigma, q), x)=\left\{
\begin{array}{ll}
m^{\lfloor\frac{q+x}{m}\rfloor}(q+_mx),&0\leq x\leq n-1,\\
m^{\lfloor\frac{q+x}{m}\rfloor},&x=n.
\end{array}
\right.$
\end{itemize}
In particular, the   values of $\psi$ are all non-empty words  if and only if $n\geq m$. Moreover, the composition $A_n\circ B_m$ is a Mealy automaton if and only if $m=n=1$.
\end{stw}
\begin{dwd}
The equalities for $X$, $Q$ and $Z$ directly follow from the definition of $A_n\circ B_m$. By Lemma~\ref{lem1}, for the transition function $\varphi$, we have:
\begin{eqnarray*}
\varphi((\sigma, q),x)&=&(\varphi^{A_n}(\sigma, x), \overline{\varphi}^{B_m}(q, \psi^{A_n}(\sigma, x)))=(\sigma, \overline{\varphi}^{B_m}(q, \psi^{A_n}(\sigma, x)))=\\
&=&\left\{
\begin{array}{ll}
(\sigma, \overline{\varphi}^{B_m}(q, 0^x1)),&0\leq x\leq n-1,\\
(\sigma, \overline{\varphi}^{B_m}(q, 0^n)),&x=n,
\end{array}
\right.=\\
&=&\left\{
\begin{array}{ll}
(\sigma, 0),&0\leq x\leq n-1,\\
(\sigma, (q+_mn)),&x=n,
\end{array}
\right.
\end{eqnarray*}
and for the output function $\psi$, we have:
\begin{eqnarray*}
\psi((\sigma, q),x)&=&\overline{\psi}^{B_m}(q, \psi^{A_n}(\sigma, x)))=\\
&=&\left\{
\begin{array}{ll}
\overline{\psi}^{B_m}(q, 0^x1),&0\leq x\leq n-1,\\
\overline{\psi}^{B_m}(q, 0^n),&x=n,
\end{array}
\right.=\\
&=&\left\{
\begin{array}{ll}
m^{\lfloor\frac{q+x}{m}\rfloor}(q+_mx),&0\leq x\leq n-1,\\
m^{\lfloor\frac{q+x}{m}\rfloor},&x=n.
\end{array}
\right.
\end{eqnarray*}

Now, if $n<m$, then $\psi((\sigma, 0),n)=m^{\lfloor n/m\rfloor}=m^0=\epsilon$. In general,  for every $0\leq q\leq m-1$ and $x\in X$ the length of the word $\psi((\sigma, q),x)$ is equal to
$$
|\psi((\sigma, q),x)|=\left\{
\begin{array}{ll}
\lfloor\frac{q+x}{m}\rfloor+1,&0\leq x\leq n-1,\\
\lfloor\frac{q+n}{m}\rfloor, &x=n.
\end{array}
\right.
$$
Hence, if $n\geq m$, then $|\psi((\sigma, q),x)|\geq 1$ for all $0\leq q\leq m-1$ and $x\in X$. If $n>m$, then
$$
|\psi((\sigma, m-1),n)|=\left\lfloor\frac{m-1+n}{m}\right\rfloor\geq 2.
$$
If $n=m>1$, then
$$
|\psi((\sigma, m-1),m-1)|=\left\lfloor\frac{m-1+m-1}{m}\right\rfloor+1\geq 2.
$$
If $n=m=1$, then $Q=\{(\sigma, 0)\}$, $X=Z=\{0,1\}$ and $\psi((\sigma,0),0)=0$, $\psi((\sigma, 0),1)=1$. The claim follows.
\end{dwd}

When  considering the composition $B_m\circ A_n$, we observe that it exits  if and only if $n=m$.

\begin{stw}\label{kl1}
For every $n\geq 1$, we have:
$$
B_n\circ A_n=(\{0,1\}, \{(0,\sigma), \ldots, (n-1, \sigma)\}, \{0,1\}, \varphi', \psi'),
$$
where
\begin{eqnarray*}
\varphi'((q, \sigma), x)&=&\left\{
\begin{array}{ll}
(q+_n1, \sigma),&x=0,\\
(0, \sigma),&x=1,
\end{array}\right.\\
\psi'((q, \sigma), x)&=&\left\{
\begin{array}{ll}
\epsilon, &x=0\;\;{\rm and}\;\;q\neq n-1,\\
0^qx, &x=1\;\;{\rm or}\;\;q=n-1.
\end{array}
\right.
\end{eqnarray*}
In particular, for all $0\leq q\leq n-1$ and $t\geq 0$, we have:
\begin{eqnarray}
\overline{\varphi'}((q, \sigma), 0^t)&=&(q+_nt,\sigma),\label{rew1}\\
\overline{\psi'}((q, \sigma), 0^t)&=&0^{n\cdot \lfloor\frac{q+t}{n}\rfloor}.\label{rew2}
\end{eqnarray}
\end{stw}
\begin{dwd}
For the transition function $\varphi'$, we can write:
\begin{eqnarray*}
\varphi'((q, \sigma), x)&=&(\varphi^{B_n}(q,x), \overline{\varphi}^{A_n}(\sigma, \psi^{B_n}(q,x)))=\\
&=&\left\{
\begin{array}{ll}
(q+_n1, \sigma),&x=0,\\
(0, \sigma),&x=1.
\end{array}\right.
\end{eqnarray*}
For the output function $\psi'$, we have:
\begin{eqnarray*}
\psi'((q, \sigma), x)&=&\overline{\psi}^{A_n}(\sigma, \psi^{B_n}(q,x)))=\\
&=&\left\{
\begin{array}{ll}
\overline{\psi}^{A_n}(\sigma, q),&x=1,\\
\overline{\psi}^{A_n}(\sigma, \epsilon), &x=0\;\;{\rm and}\;\;q\neq n-1,\\
\overline{\psi}^{A_n}(\sigma, n), &x=0\;\;{\rm and}\;\;q=n-1,
\end{array}
\right.=\\
&=&\left\{
\begin{array}{ll}
0^q1,&x=1,\\
\epsilon, &x=0\;\;{\rm and}\;\;q\neq n-1,\\
0^n, &x=0\;\;{\rm and}\;\;q=n-1,
\end{array}
\right.=\\
&=&\left\{
\begin{array}{ll}
\epsilon, &x=0\;\;{\rm and}\;\;q\neq n-1,\\
0^qx, &x=1\;\;{\rm or}\;\;q=n-1.
\end{array}
\right.
\end{eqnarray*}

The equality (\ref{rew1}) directly follows from the above formula for $\varphi'$. To show (\ref{rew2}), we can write
$$
\overline{\psi'}((q, \sigma), 0^t)=x_1x_2\ldots x_t,
$$
where
$$
x_i=\psi'((q+_n(i-1), \sigma),0),\;\;\;i=1,2,\ldots, t.
$$
Hence, by the above formula for $\psi'$, we have: $x_i\in\{\epsilon, 0^n\}$, and $x_i=0^n$ if and only if $q+_n(i-1)=n-1$.
Thus
$$
\overline{\psi'}((q,\sigma), 0^t)=0^{nr},
$$
where $r\geq 0$ is the number of those elements from the sequence
$$
q,\;q+_n1,\;q+_n2,\;\ldots, q+_n(t-1)
$$
which are equal to $n-1$. The similar argument as in the proof of Lemma~\ref{lem1} gives the equality $r=\lfloor\frac{q+t}{n}\rfloor$. The claim follows.
\end{dwd}

\begin{theorem}\label{twa3}
If $n\geq m$, then for every $0\leq q\leq m-1$ the function
$$
f_{\omega, (\sigma, q)}^{A_n\circ B_m}\colon \{0,1,\ldots, n\}^\omega\to \{0,1,\ldots, m\}^\omega
$$
is defined by a Mealy automaton, i.e. $f_{\omega, (\sigma, q)}^{A_n\circ B_m}$ is a short map. Moreover, for all $n,m\geq 1$ the function $f_{\omega, (\sigma, 0)}^{A_n\circ B_m}$ is bijective.
\end{theorem}
\begin{dwd}
The first part of the claim directly follows from Proposition~\ref{ppp1} and Corollary~\ref{cor22}. To show the second part, we use the previous observation that both the automata $A_n$ and $B_m$ are nondegenerate. Now, by Proposition~\ref{ppp2},  it is enough to show that for every $n\geq 1$ both the functions
$$
f_{\omega, \sigma}^{A_n}\colon \{0,1,\ldots, n\}^\omega\to\{0,1\}^\omega
$$
and
$$
f_{\omega, 0}^{B_n}\colon \{0,1\}^\omega\to\{0,1,\ldots, n\}^\omega
$$
are bijective. It can be verified by just looking at the diagrams of the automata $A_n$ and $B_n$. Our formal proof below additionally shows that the functions $f_{\omega, \sigma}^{A_n}$ and $f_{\omega, 0}^{B_n}$ are mutually inverse, i.e. the following two equalities hold:
$$
f_{\omega, (\sigma, 0)}^{A_n\circ B_n}=Id_{\{0,1,\ldots, n\}^\omega},\;\;\;f_{\omega, (0, \sigma)}^{B_n\circ A_n}=Id_{\{0,1\}^\omega}.
$$
To show the first equality, it is enough to show that the automaton $A_n\circ B_n$ being in the initial state $(\sigma, 0)$ and reading any letter $x\in \{0,1,\ldots, n\}$ from the input tape, it remains at the state $(\sigma, 0)$ and sends the same letter $x$ on the output tape. But this easily follows from the formulae for $\varphi$ and $\psi$ in Proposition~\ref{ppp1}. Indeed, we have by these formulae:
$$
\varphi((\sigma, 0), x)=\left\{
\begin{array}{ll}
(\sigma, 0),&0\leq x\leq n-1,\\
(\sigma, 0+_nn),&x=n,
\end{array}\right.
$$
and
$$
\psi((\sigma, 0), x)=\left\{
\begin{array}{ll}
n^{\lfloor\frac{0+x}{n}\rfloor}(0+_nx),&0\leq x\leq n-1,\\
n^{\lfloor\frac{0+x}{n}\rfloor},&x=n.
\end{array}
\right.
$$
Thus $\varphi((\sigma, 0), x)=(\sigma, 0)$ and $\psi((\sigma, 0), x)=x$ for any $x\in\{0,1,\ldots, n\}$.

To show the equality $f_{\omega, (0, \sigma)}^{B_n\circ A_n}=Id_{\{0,1\}^\omega}$, we use the fact that every infinite word $u\in\{0,1\}^\omega$ in the metric  space $(\{0,1\}^\omega, d_\lambda)$ is the limit of the sequence $(v_i)_{i\geq 1}$, where $v_i=u_i1^\infty$ and $u_i$ is the prefix of length $|u_i|=i $ of the word $u$. Hence,  since the function $f_{\omega, (0, \sigma)}^{B_n\circ A_n}$ is continuous, it is enough to show the equalities: $f_{\omega, (0, \sigma)}^{B_n\circ A_n}(u_i1^\infty)=u_i1^\infty$ for every $i\geq 1$. For this aim, we show  that $f_{(0,\sigma)}^{B_n\circ A_n}(w1)=w1$ for every finite word $w\in\{0,1\}^*$.
Since  every finite word over $\{0,1\}$ with the last letter $1$ is a concatenation of finitely many  words of the form $0^t1$ ($t\geq 0$),  it is enough to show that for every $t\geq 0$ the automaton $B_n\circ A_n$ works as follows: being in the initial state $(0,\sigma)$ and reading the word $0^t1$ from the input tape, it returns  to the initial state $(0, \sigma)$ and sends the same word $0^t1$ on the output tape. In other words, it is enough to show the equalities:
\begin{eqnarray*}
\overline{\varphi'}((0, \sigma), 0^t1)&=&(0,\sigma),\\
\overline{\psi'}((0, \sigma), 0^t1)&=&0^t1.
\end{eqnarray*}
But the above two equalities easily follow from the formulae for $\varphi'$ and $\psi'$ and from the equalities~(\ref{rew1})--(\ref{rew2}) in Proposition~\ref{kl1}. Indeed, since $\varphi'((q,\sigma),1)=(0,\sigma)$ for every $0\leq q\leq n-1$, we obtain by (\ref{rew1}):
\begin{eqnarray*}
\overline{\varphi'}((0, \sigma), 0^t1)=\varphi'(\overline{\varphi'}((0,\sigma), 0^t),1)=\varphi'((0+_nt,\sigma),1)=(0,\sigma).
\end{eqnarray*}
Next, if we denote by $[t]_n\in\{0,1,\ldots, n-1\}$ the remaining from the division of $t$ by $n$, then the formulae~(\ref{rew1})--(\ref{rew2}) and the equality $\psi'(([t]_n,\sigma),1)=0^{[t]_n}1$ imply:
\begin{eqnarray*}
\overline{\psi'}((0, \sigma), 0^t1)&=&\overline{\psi'}((0, \sigma),0^t)\psi'(\overline{\varphi'}((0,\sigma), 0^t),1)=\\
&=&0^{n\cdot \lfloor\frac{0+t}{n}\rfloor}\psi'(([t]_n, \sigma),1)=\\
&=&0^{n\cdot \lfloor\frac{0+t}{n}\rfloor}0^{[t]_n}1=0^{n\cdot\lfloor\frac{t}{n}\rfloor+[t]_n}1=0^t1.
\end{eqnarray*}
The claim follows.
\end{dwd}

\begin{figure}[hbtp]
\centering
\includegraphics[width=13cm]{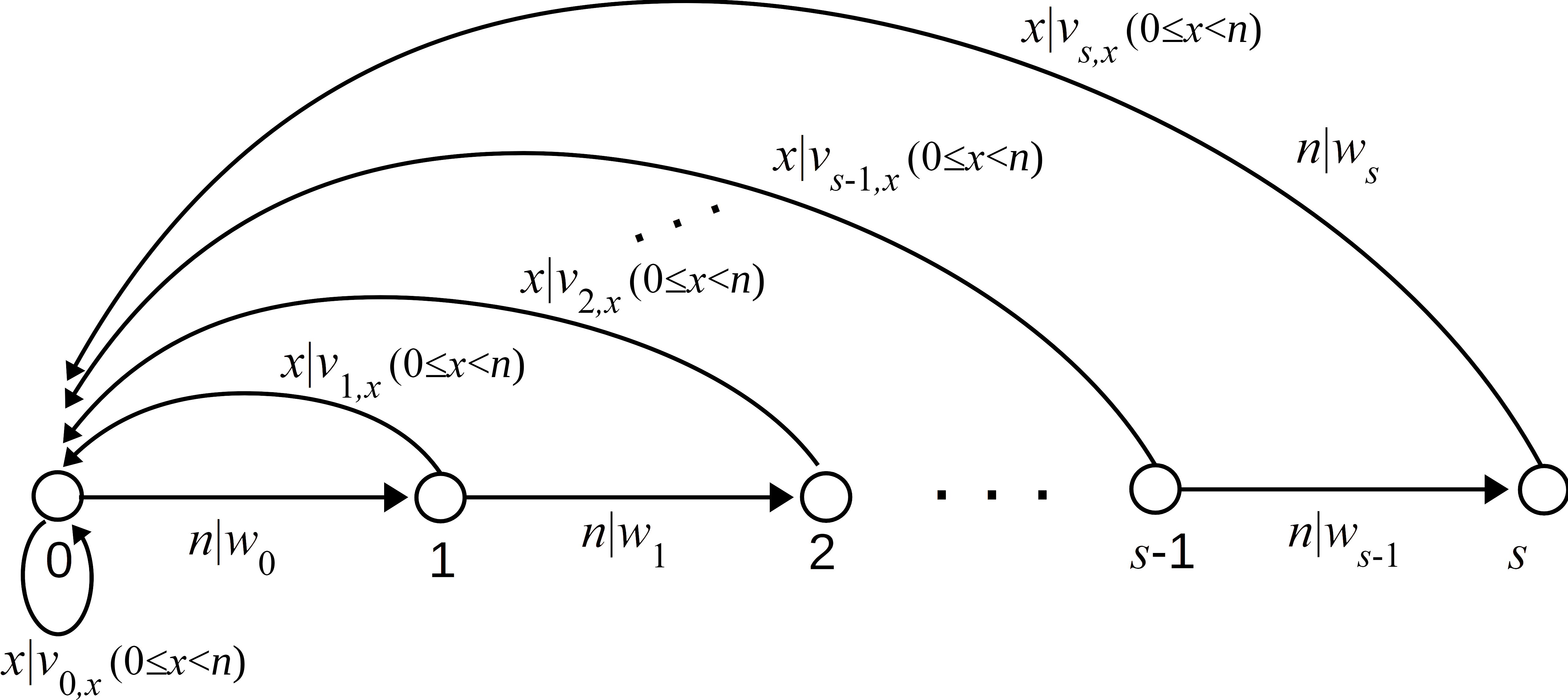}
\caption{The automaton $C$: $s=\frac{m}{{\rm gcd}(m,n)}-1$,\;\;\;$v_{q,x}=m^{\left\lfloor\frac{[qn]_m+x}{m}\right\rfloor}(qn+_mx)$,\;\;\; $w_q=m^{\left\lfloor\frac{[qn]_m+n}{m}\right\rfloor}$}
\label{fig3}
\end{figure}

If we take a closer look at the transition function $\varphi$ of the automaton $A_n\circ B_m$ from Proposition~\ref{ppp1}, then we see that the subset
$$
Q'=\left\{(\sigma, [in]_m)\colon 0\leq i<\frac{m}{{\rm gcd}(n,m)}\right\}\subseteq Q=\{(\sigma, 0), (\sigma, 1), \ldots, (\sigma, m-1)\}
$$
induces a subautomaton (connected component) of the automaton $A_n\circ B_m$, i.e. $\varphi(q,x)\in Q'$ for all $q\in Q'$, $x\in X$. As a result, we obtain the following construction of the automaton $C$ defining in the case $n\geq m$ a bijective short map between the Cantor spaces
$\{0,1\ldots, n\}^\omega$ and $\{0,1,\ldots, m\}^\omega$. The Moore'a diagram of the automaton $C$ is depicted in Figure~\ref{fig3}.

\begin{wn}\label{wn22}
Let $m,n\geq 1$,  $\eta_{m,n}:=\frac{m}{{\rm gcd}(n,m)}$, and let
$$
C=(\{0,1,\ldots, n\}, Q^C, \{0,1,\ldots, m\},\varphi^C, \psi^C)
$$
be an asynchronous automaton in which
\begin{itemize}
\item $Q^C=\left\{0,1,\ldots, \eta_{m,n}-1\right\}$,
\item $\varphi^C(q, x)=\left\{
\begin{array}{ll}
0,&0\leq x\leq n-1,\\
q+_{\eta_{m,n}}1,&x=n,
\end{array}\right.$
\item $\psi^C(q, x)=\left\{
\begin{array}{ll}
m^{\left\lfloor\frac{[qn]_m+x}{m}\right\rfloor}([qn]_m+_mx),&0\leq x\leq n-1,\\
m^{\left\lfloor\frac{[qn]_m+x}{m}\right\rfloor},&x=n.
\end{array}
\right.$
\end{itemize}
Then $f_{\omega, 0}^C=f_{\omega, (\sigma, 0)}^{A_n\circ B_m}$. In particular, if $n\geq m$, then  $f_{\omega, 0}^C$ is a bijective short map between the Cantor spaces
$\{0,1\ldots, n\}^\omega$ and $\{0,1,\ldots, m\}^\omega$.
\end{wn}

\section{Finite Mealy automata defining bijective colouring}\label{sec5}

The construction of the bijective  function
$$
f\colon \{0,1,\ldots, n\}^\omega\to \{0,1,\ldots, m\}^\omega
$$
defined by a Mealy automaton from the previous section is based on  the asynchronous automaton $A_n\circ B_m$ with the $m$-element set of states. According to Corollary~\ref{wn22}, the automaton $A_n\circ B_m$ contains a subautomaton with $\eta_{m,n}=m/{\rm gcd}(m,n)$ states, which also defines $f$.  In particular, if $m$ divides $n$, then this subautomaton is a 1-state asynchronous automaton. In this context, it is natural to ask if there is a finite Mealy automaton $A=(X, Q, Y, \varphi, \psi)$ with $|X|\neq |Y|$ and such  that the function  $f_{\omega, q}^A\colon X^\omega\to Y^\omega$ is bijective. As we previously observed, if $|X|<|Y|$, then  $f_{\omega, q}^A$ cannot be surjective. We prove now the  following statement.

\begin{theorem}\label{stst1}
Let $A=(X, Q, Y, \varphi, \psi)$ be an arbitrary Mealy automaton such that $|X|>|Y|$ and the function  $f_{\omega, q}^A\colon X^\omega\to Y^\omega$ is injective for some state $q\in Q$. Then $|Q|=\infty$.
\end{theorem}
\begin{dwd}
Let us denote $f=f_{\omega, q}^A$. The following lemmas hold:
\begin{lem}\label{lemlem1}
For every  $k\geq 0$ there are  $w\in X^*$ and $u,u'\in X^\omega$ such that $P(u,u')=\epsilon$ and $|P(f(wu),f(wu'))|\geq |w|+k$.
\end{lem}
\begin{dwd}(of Lemma~\ref{lemlem1})
Let $k\geq 0$ be arbitrary. Since $|X|>|Y|$, there is $n\geq 0$ such that $|X|^{n}>|Y|^{n+k}$. Let us fix $v_0\in X^k$. Obviously, for every  $v\in X^n$ and $u\in X^\omega$ the prefix of length $n+k$ of the word $f(vv_0u)$  belongs to the set  $Y^{n+k}$. Since $|X^n|>|Y^{n+k}|$, there are   $v,v'\in X^n$ and $u_0,u_1\in X^\omega$ such that $v\neq v'$ and the words  $f(vv_0u_0)$ and $f(v'v_0u_1)$ have a common prefix of length  $n+k$, i.e.
\begin{equation}\label{aj11}
|P(f(vv_0u_0), f(v'v_0u_1))|\geq n+k.
\end{equation}
Let $w:=P(v,v')$. Since  $v\neq v'$ and $|v|=|v'|=n$, we obtain: $|w|<n$. Hence  $v=ww_0$, $v'=ww_1$ for some $w_0, w_1\in X^*\setminus\{\epsilon\}$. In particular, if we denote  $u:=w_0v_0u_0$ and $u':=w_1v_0u_1$, then the following equalities hold:
$$
vv_0u_0=ww_0v_0u_0=wu,\;\;\;v'v_0u_1=ww_1v_0u_1=wu'.
$$
Thus, by (\ref{aj11}), we can write: $|P(f(wu),f(wu'))|\geq n+k>|w|+k$. Since the words $w_0$, $w_1$ are non-empty and $P(w_0, w_1)=\epsilon$, we also have: $P(u,u')=\epsilon$.
\end{dwd}

\begin{lem}\label{lem2}
For every  $w\in X^*$ there is  $k\geq 0$ such that for all $u,u'\in X^\omega$ the equality $P(u,u')=\epsilon$ implies $|P(f(wu),f(wu'))|<|w|+k$.
\end{lem}
\begin{dwd}(of Lemma~\ref{lem2})
Suppose, contrary, that there is $w\in X^*$ such that for every  $k\geq 0$ there are $u_k, u'_k\in X^\omega$ which satisfy:
$$
P(u_k, u'_k)=\epsilon\;\;\;{\rm and}\;\;\;|P(f(wu_k), f(wu'_k))|\geq |w|+k.
$$
Since the space  $(X^\omega, d_{X,\lambda})$ is compact, the sequence  $(u_k)_{k\geq 0}$ contains an infinite convergent subsequence $(u_{k_i})_{i\geq 0}$:
$$
u_{k_i}\to v,\;\;\;v\in X^\omega.
$$
Obviously, the corresponding subsequence $(u'_{k_i})_{i\geq 0}$ of the sequence  $(u'_k)_{k\geq 0}$ need not be convergent, but it contains some infinite convergent subsequence $(u'_{l_i})_{i\geq 0}$:
$$
u'_{l_i}\to v',\;\;\;v'\in X^\omega.
$$
Clearly, the sequence $(u_{l_i})_{i\geq 0}$ is an infinite subsequence of the sequence $(u_{k_i})_{i\geq 0}$, and hence:
$$
u_{l_i}\to v.
$$
Since $P(u_{l_i}, u'_{l_i})=\epsilon$ for every $i\geq 0$, we obtain: $P(v,v')=\epsilon$, which implies $v\neq v'$. Since
$$
wu_{l_i}\to wv,\;\;\;wu'_{l_i}\to wv',
$$
and $f$ is continuous, we have:
\begin{equation}\label{c8}
f(wu_{l_i})\to f(wv),\;\;\;f(wu'_{l_i})\to f(wv').
\end{equation}
For every  $i\geq 0$, we also have: $|P(f(wu_{l_i}), f(wu'_{l_i}))|\geq |w|+l_i$, and hence:
$$
|P(f(wu_{l_i}), f(wu'_{l_i}))|\to\infty,
$$
and by the two previous  convergences, we obtain: $f(wv')=f(wv)$. Since $f$ is injective, we have: $wv=wv'$, which implies $v=v'$, contrary to our previous observation.
\end{dwd}

By Lemmas~\ref{lemlem1}--\ref{lem2}, there is an infinite strictly increasing sequence $(k_i)_{i\geq 0}$ of numbers and an infinite sequence of words $(w_i)_{i\geq 0}$  which satisfy: \begin{itemize}
\item $k_0=0$,
\item for every $i\geq0$ the word $w_i\in X^*$ corresponds to the number $k:=k_i$ from Lemma~\ref{lemlem1},
\item for every $i\geq 0$ the number $k_{i+1}$ is any number greater than $k_i$ which corresponds to the word $w:=w_i$ from Lemma~\ref{lem2} (note that if a number $k$ satisfies Lemma~\ref{lem2}, then  every number greater than $k$ also satisfies this lemma).
\end{itemize}
In particular, there are sequences $(u_i)_{i\geq 0}$ and $(u'_i)_{i\geq 0}$ of infinite words $u_i, u'_i\in X^\omega$ such that $P(u_i,u'_i)=\epsilon$ for every $i\geq 0$, and for all $i,j\geq 0$ the inequalities hold:
\begin{eqnarray*}
|P(f(w_iu_i),f(w_iu'_i))|&\geq& |w_i|+k_i,\\
|P(f(w_iu_j),f(w_iu'_j))|&<&|w_i|+k_{i+1}.
\end{eqnarray*}

\begin{lem}\label{lem3}
If $i\neq j$, then $f_{\omega, q_i}^A\neq f_{\omega, q_j}^A$, where $q_i=\overline{\varphi}(q, w_i)\in Q$ for every $i\geq 0$.
\end{lem}
\begin{dwd}(of Lemma~\ref{lem3})
We can assume $i<j$. Let us denote
$$
u:=u_j,\;\;u':=u'_j,\;\;g:=f_{\omega, q_i}^A,\;\;h:=f_{\omega, q_j}^A.
$$
In particular, the following inequalities hold
\begin{eqnarray}
|P(f(w_ju), f(w_ju'))|&\geq&|w_j|+k_j,\label{k12}\\
|P(f(w_iu), f(w_iu'))|&<&|w_i|+k_{i+1}.\label{k13}
\end{eqnarray}
We also have:
\begin{eqnarray*}
f(w_ju)=f_q^A(w_j)h(u),\;\;\;f(w_ju')=f_q^A(w_j)h(u'),\\
f(w_iu)=f_q^A(w_i)g(u),\;\;\;f(w_iu')=f_q^A(w_i)g(u').
\end{eqnarray*}
Since $|f_q^A(w_j)|=|w_j|$ and $|f_q^A(w_i)|=|w_i|$, we obtain by (\ref{k12})--(\ref{k13}):
$$
|P(h(u), h(u'))|\geq k_j,\;\;|P(g(u), g(u'))|< k_{i+1}.
$$
Since $k_j\geq k_{i+1}$, we can write:
$$
|P(h(u), h(u'))|>|P(g(u), g(u'))|,
$$
which implies $g\neq h$.
\end{dwd}

By Lemma~\ref{lem3}, we obtain that the subset $\{\overline{\varphi}(q,w_i)\colon i\geq 0\}\subseteq Q$ is infinite. Hence $|Q|=\infty$, which finishes the proof of Theorem~\ref{stst1}.
\end{dwd}

\begin{wn}\label{wnwn1}
Let $A=(X, Q, Y, \varphi, \psi)$ be a finite Mealy automaton such that the function $f_{\omega, q}^A\colon X^\omega\to Y^\omega$ is bijective for some state $q\in Q$. Then $|X|=|Y|$. In particular, if $|X|\neq |Y|$, then there is no bijective short map $f\colon X^\omega\to Y^\omega$ defined by a finite Mealy automaton.
\end{wn}

\end{document}